\documentclass[twoside]{article}
\begin{document}
\vspace*{-2ex}
\noindent{\footnotesize {\em Sankhy\=a\,: The Indian Journal of 
   Statistics}\\[-1ex]
1999, Volume 61, Series A, Pt. xx, pp. xx-xx} \\ \vspace{2ex}

\markboth {\hfill \sc m.c. jones, d.f. signorini and n.l. hjort \hfill}
    {\hfill \sc multiplicative bias correction in kernel density estimation \hfill}

\begin{center}
{\large ON MULTIPLICATIVE BIAS CORRECTION IN KERNEL DENSITY ESTIMATION} \\
\vspace{1em}
{\it By} \ M.C. JONES\\
{\it The Open University, U.K.} \\
D.F. SIGNORINI\\
{\it Quintiles Scotland Ltd, U.K.} \\
and \\
N.L. HJORT\\
{\it University of Oslo, Norway} \\
\end{center}
\vspace{1em}

{\footnotesize {\it SUMMARY.} 
Hjort and Glad (1995) present a method for semiparametric density estimation.
Relative to the ordinary kernel density estimator, this technique performs
much better when a parametric vehicle distribution fits the data, and 
otherwise
performs at broadly the same level. Jones, Linton and Nielsen (1995) present
a somewhat similar method for density estimation which has higher order bias 
for all sufficiently smooth densities. In this paper, we combine the two
methods. We show that, theoretically, the desired properties of general 
higher order bias allied with even better performance for an appropriate 
vehicle model are achieved. Simulations suggest that the new estimator
realises only a little of its theoretical potential in practice for small 
to moderately large sample sizes.}

\vspace{1em}
\begin{center}
1. \quad {\bf Introduction} 
\end{center}
\vspace{1em}

Two promising recent proposals for `improved' kernel density estimation share a
common form but exhibit rather different types of performance. In this paper,
we investigate combining the two approaches in an attempt to obtain the
best of both worlds.

The common formulation is as follows. Introduce a
kernel function $K$ which we will take to be a symmetric
probability density function, and its associated smoothing parameter,
or bandwidth, $h$, writing $K_h(u) \equiv h^{-1}K(h^{-1}u)$.
Let $g$ be a function to be specified.
Then, for a random sample $X_1, \ldots ,X_n$, of size $n$,
consider estimators of the density $f$ of the form
$$\tilde{f}(x) = g(x) n^{-1} \sum_{i=1}^n g(X_i)^{-1}K_h(x-X_i) \eqno \ldots (1.1)$$

\ \hrule 
\vspace{.5em}
{\footnotesize \noindent Paper received. June 1998; revised June 1999. \\
{\it AMS} (1991) {\it subject classification}: 62G07. \\
{\it Key words and phrases.} Bias reduction; semiparametric density 
estimation; smoothing.}

\newpage

\noindent
and also its renormalisation to achieve unit integral
$${\tilde{f}(x) \over \int\tilde{f}(z)dz} = {g(x)  \sum_{i=1}^n
g(X_i)^{-1}K_h(x-X_i) \over \sum_{i=1}^n
g(X_i)^{-1} (K_h*g)(X_i)} \eqno \ldots (1.2)$$
where $*$ denotes convolution.

The simplest special case is when $g(x) \equiv 1$, in which case we get the
basic kernel density estimator (KDE)
$$\hat{f}(x) = n^{-1} \sum_{i=1}^n K_h(x-X_i)  \eqno \ldots (1.3)$$
(which automatically integrates to one).
Theoretical properties of this estimator are well known (e.g.\ Scott, 1992,
Wand and Jones, 1995). In particular, provided that $f$ has two continuous
derivatives, as $n \rightarrow \infty$ and $h = h(n)
\rightarrow 0$, the bias of $\hat{f}$ is of order $h^2$ and its variance is
$O((nh)^{-1})$ (provided also that $nh \rightarrow \infty$).

If $g$ is taken to be an initial, or `pilot', estimator of $f$, then $\tilde{f}$
becomes a two-stage multiplicatively corrected density estimator. In
particular, such an estimator acts as a multiplicative bias correction:
the two appearances of $g$ in (1.1) are such that the leading bias in $g$
as an estimator of $f$ occurs with opposite signs and cancels out.

If we take $g(x) = \hat{f}(x)$, we obtain the estimator of
Jones, Linton and Nielsen (1995) i.e.\
$$\hat{f}_N(x) = \hat{f}(x) n^{-1} \sum_{i=1}^n \hat{f}(X_i)^{-1}
K_h(x-X_i)  \eqno \ldots (1.4)$$
or in renormalised form
$$\hat{f}_N^R(x) = {\hat{f}(x) \sum_{i=1}^n \hat{f}(X_i)^{-1}
K_h(x-X_i) \over \sum_{i=1}^n \hat{f}(X_i)^{-1} 
(K_h*\hat{f})(X_i) } . \eqno \ldots (1.5)$$
The subscript $N$ stands for (fully) Nonparametric and the superscript $R$,
when present, denotes Renormalisation.

The bias cancellation here works to afford a bias of order $h^4$, provided
we now assume that $f$ has four continuous derivatives. We thus refer
to $\hat{f}_N$ as a higher order bias kernel density estimator or HOBKDE.
Asymptotic variance of $\hat{f}_N$ remains of order $(nh)^{-1}$ although
there is an increase in the constant coefficient. The achievement of decreased
bias at the expense, in finite samples,
of increased variance is thus somewhat disguised. A great many HOBKDE proposals
have been made; for a review and comparison see Jones and Signorini (1997).
The evidence is that this idea, at least in its $\hat{f}_N^R$ form, is amongst
the best HOBKDEs.

On the other hand, let $f(x; \hat{\theta})$ be one of the usual parametric fits
of the parametric family $f(x;\theta)$ to the data. Hjort and Glad (1995)
proposed (1.1) using $g(x) = f(x;\hat{\theta})$ i.e.\
$$\hat{f}_S(x) = f(x;\hat{\theta}) n^{-1} \sum_{i=1}^n f(X_i;\hat{\theta})^{-1}
K_h(x-X_i) . \eqno \ldots (1.6)$$
Its renormalised form is, of course,
$$\hat{f}_S^R(x) = {f(x;\hat{\theta}) \sum_{i=1}^n f(X_i;\hat{\theta})^{-1}
K_h(x-X_i) \over \sum_{i=1}^n f(X_i;\hat{\theta})^{-1}
(K_h* f(\cdot;\hat{\theta}))(X_i) } . \eqno \ldots (1.7)$$

In general, the bias in using $f(x;\hat{\theta})$ is $O(1)$, being so unless
$f$ happens to belong to the parametric class $f(\cdot;\theta)$. The bias
correction works to cancel the $O(1)$ biases with the end result that
the bias in using $\hat{f}_S$ (or $\hat{f}_S^R$)
is of order $h^2$. However, if the parametric model does encompass $f$,
the $O(h^2)$ bias term also vanishes; in fact, the multiplier of $h^2$ depends
on $(f/f_0)''(x)$ where $f_0(x)$ is that version of $f(x;\theta)$ 
which is closest to the true density $f$ in a sense appropriate
to the particular parametric method being used. (Moreover,
the bias is greatly reduced in this case, all that remains being any bias due
to estimating $\theta$ by $\hat{\theta}$.) The asymptotic variance of
$\hat{f}_S$ and $\hat{f}_S^R$ remains precisely the same, constants included,
as that of $\hat{f}$.

The subscript S attached to the Hjort and Glad (1995) estimator 
stands for Semiparametric, and $\hat{f}_S$ and $\hat{f}_S^R$ are
examples of semiparametric KDEs (SKDEs). SKDEs attempt to obtain the best
aspects of both parametric and nonparametric density estimation: when one
has a good parametric model for the data, the aim is to achieve
the greater efficiency of parametric model fitting; if one's parametric
proposal proves not to be a good one, the method `becomes nonparametric'
and should still perform well --- around the level of $\hat{f}$
--- whatever is $f$.
The bandwidth plays a major role here: large $h$ corresponds essentially
(provided we renormalise) to fitting the parametric model, but when $h$ is
small, the kernel smoothing side takes over.

Also, $\hat{f}_S$ and $\hat{f}_S^R$ are part of a plethora of SKDE proposals,
currently being reviewed and compared by Hjort, Jones and Storvik (paper in
preparation). Early indications are that $\hat{f}_S$ and $\hat{f}_S^R$ are
among the best such methods.

Can we obtain both semiparametric performance, in particular leading
bias zeroed for a parametric family, and HOBKDE performance, that is,
bias at arbitrary $f$ of $O(h^4)$ (sufficient smoothness of $f$ permitting)?
The answer is affirmative, and we exhibit one way of accomplishing this by
combining the two methods described above. The new method is basically to
set $g = \hat{f}_S$ in (1.1); see (2.1).

We explore the asymptotic bias and variance properties of our proposal, and
of variations thereon incorporating renormalisation, in Section 2. 
Simulations comparing the new methods with the ordinary KDE and with the
Jones, Linton and Nielsen and Hjort and Glad estimators 
are made in Section 3.
In the simulations, a normal parametric model only is used along with a 
practically unavailable optimal bandwidth selection. The end result, from 
the practical viewpoint, is a little disappointing in this case.
For $n=100$, it is $\hat{f}_N^R$ that dominates, and not the new higher
order bias semiparametric density estimator;
%$\hat{f}_{S,N}$;
for $n=500$, the two have similar performance. We 
make our brief conclusions in Section 4. 

\vspace{1em}
\begin{center}
2. \quad {\bf Method and Theoretical Results}
\end{center}
\vspace{1em}

The novel higher order bias semiparametric kernel density estimator 
(HOBSKDE) proposed
in this paper results from combining the Jones, Linton and Nielsen (1995)
higher order bias and Hjort and Glad (1995) semiparametric density
estimators in the following way. Define
$$\hat{f}_{S,N} (x) =  \hat{f}_S(x) n^{-1} \sum_{i=1}^n \hat{f}_S(X_i)^{-1}
K_h(x-X_i)  \eqno \ldots (2.1)$$
and its renormalised form
$$\hat{f}_{S,N}^R(x) = {\hat{f}_S(x) \sum_{i=1}^n \hat{f}_S(X_i)^{-1}
K_h(x-X_i) \over \sum_{i=1}^n \hat{f}_S(X_i)^{-1} (K_h
*\hat{f}_S)(X_i) } . \eqno \ldots (2.2)$$
(Note that renormalisation of $\hat{f}_S(x)$ makes no difference because
the renormalisation constant cancels out.)

For the general form (1.1), the mean is easily seen to be
\begin{eqnarray}
g(x) \{K_h*(f/g)\}(x)
 &\simeq& g(x) \left\{\left({f \over g}\right)(x) + \frac{h^2}2 s_2
\left( {f \over g} \right)''(x) +\frac{h^4}{24} s_4 
\left( {f \over g} \right)^{(iv)}(x) \right\} \nonumber \\
&=& f(x) + \frac{h^2}2 s_2 g(x) \left({f \over g}\right)''(x) + 
\frac{h^4}{24} s_4 g(x) \left( {f \over g} \right)^{(iv)}(x)  \ldots (2.3) 
\nonumber \end{eqnarray}
as $h \rightarrow 0$. Here, $s_\ell = \int u^\ell K(u)du$. 
Ignoring the difference between $\hat{\theta}$ and $\theta$ which is
negligible for these asymptotic purposes (using a variation of arguments 
used in Hjort and Glad, 1995, Section 3), a special case of this is that
$$E\{\hat{f}_S(x)\}- f(x) \simeq \frac{h^2}2 s_2 f_0(x) (f/f_0)''(x).
\eqno \ldots (2.4)$$
Another special case arises by replacing $g(x)$ in (2.3) by the expansion
for $E\{\hat{f}(x)\}$
to give the bias expression for (1.4) given by Jones, Linton and Nielsen
(1995). There are slightly different expressions for renormalised forms 
which will not be given here.

Inserting (2.4) in (2.3) --- which turns out to give the same answer as 
a more rigorous calculation --- yields 
\begin{eqnarray}
 E\{\hat{f}_{S,N}(x)\}-f(x)& \simeq&
\frac{h^2}2 s_2 f(x) \left\{1-\frac{h^2}2 s_2 \left({f_0 \over f} 
\right)(x) \left({f \over f_0}\right)''(x)\right\}'' 
%+ \frac{h^4}{24} s_4 f(x)1^{(iv)}  
\nonumber \\ 
&=& -\frac{h^4}{4} s_2^2 f(x)\left\{\left({f_0 \over f}\right)(x) 
\left({f \over f_0}\right)''(x)\right\}'' ~~~~~~~~~~~~~~~\ldots (2.5)\nonumber 
\end{eqnarray}
Notice that the properties of a HOBSKDE pertain to $\hat{f}_{S,N}$ by (2.5):
it has $O(h^4)$ bias whatever is $f$ (sufficient smoothness permitting), and
the leading bias is zeroed when the `right' parametric family is chosen i.e.\
$f_0 = f$.

It is not difficult to see that renormalisation leads to
\begin{eqnarray}
 E\{\hat{f}_{S,N}^R(x)\}-f(x) & \simeq &
 -\frac{h^4}{4} s_2^2 f(x) \left[ \left\{\left({f_0 \over f}\right)(x) 
\left({f \over f_0}\right)''(x)\right\}'' \right. \nonumber \\
&~&~~~~~- \int \left.
 f(z)\left\{\left({f_0 \over f}\right)(z) \left({f 
\over f_0}\right)''(z)\right\}'' dz 
\right] . \nonumber
\end{eqnarray}

A careful calculation parallel to that on pp. 337-8 of Jones, Linton and
Nielsen (1995) gives the asymptotic variance of $\hat{f}_{S,N}$. It is
$$V\{\hat{f}_{S,N}(x)\} \simeq (nh)^{-1} f(x) \int (2K(u)-K*K(u))^2 du .
\eqno \ldots (2.6) $$
This corresponds exactly to the asymptotic variance of $\hat{f}_N(x)$.
(The same applies to $V\{\hat{f}^R_{S,N}(x)\}$.) One can interpret this
as $\hat{f}_{S,N}$ exhibiting a semiparametric yet $O(h^4)$ bias while
retaining the variance of the nonparametric $O(h^4)$ bias method $\hat{f}_N$.

The mean squared error (MSE) of $\hat{f}_{S,N}$ follows by adding the 
square of
(2.5) to (2.6). As with nonparametric higher order bias methods, the best
achievable rate of convergence of MSE is $O(n^{-8/9})$ in general, when
$h \sim n^{-1/9}$, but with better performance when the parametric vehicle
model is correct. 

We have not attempted to develop automatic bandwidth selection based on
these results, although they give a clear potential for `plug-in' 
estimation. There are two main reasons for this, each relegating the 
problem to a position down our list of priorities. The first is that the 
(simpler) equivalent problem for SKDEs has not yet been addressed in
detail. The second is that the relative negativeness of the results to 
come question the importance of pursuing this course. We point out,
however, that the arguments used to justify the cross validation method 
for ordinary kernel density estimation also apply to the new estimator.

\newpage
\vspace{1em}
\begin{center}
{\bf 3. Simulation Results}
\end{center}
\vspace{1em}

We follow Jones and Signorini (1997) in providing practical comparisons based
on simulations from a set of ten known densities. These densities are the
first ten normal mixtures in Figure 1 of Marron and Wand (1992). They are
referred to as ``Gaussian'', ``Skewed Unimodal'', ``Strongly Skewed'',
``Kurtotic Unimodal'', ``Outlier'', ``Bimodal'', ``Separated Bimodal'', ``Skewed
Bimodal'', ``Trimodal'' and ``Claw'', respectively.
One thousand random samples of sizes $n=100$ and $n=500$
were generated from each distribution.

\vspace{1em}

{\footnotesize

\begin{center}

\begin{tabular}{p{3.4in}}
Table 1. {\sc Means and standard errors of minimised ISE $\times 10^5$,
for samples of size $n=100$ and $n=500$ from each of the first ten 
Marron--Wand densities, over 1000 simulations for each of the following 
estimators: $\hat{f}$, the basic kernel density estimator (1.3);
$\hat{f}_N^R$, the renormalised Jones, Linton and Nielsen estimator (1.5);
$\hat{f}_S$, the raw Hjort and Glad estimator (1.6);
$\hat{f}_{S,N}$, the raw HOBSKDE estimator (2.1);
$\hat{f}_{S,N}^R$, the renormalised HOBSKDE estimator (2.2).
}
\end{tabular}

\vspace{1em}

\begin{tabular}{c|c|c|c}
\hline
 & & $n=100$ & $n=500$ \\
Density & Estimator & Mean Min.\ ISE  &
Mean Min.\ ISE  \\
 & & (S.E.)  & (S.E.)   \\
\hline
~ & $\hat{f}$ & 462~(12)  & 154~(3)   \\
    & $\hat{f}_N^R$ & 219~(~7)&  ~58~(2)  \\
Gaussian & $\hat{f}_S$ & 226~(~7)  & ~47~(1) \\
& $\hat{f}_{S,N}$ & 293~(~7)& ~74~(2) \\
 & $\hat{f}_{S,N}^R$ & 263~(~7) & ~60~(1)  \\
\hline
~ & $\hat{f}$ & 755~(17) & 234~(5)  \\
Skewed & $\hat{f}_N^R$ & 477~(13)& 135~(4) \\
Unimodal &  $\hat{f}_S$ & 605~(14)  &176~(4) \\
& $\hat{f}_{S,N}$ & 584~(13)& 141~(3) \\
 & $\hat{f}_{S,N}^R$ & 574~(13) & 134~(3)  \\
\hline
~ & $\hat{f}$ & 4227~(53) & 1345~(15)  \\
Strongly& $\hat{f}_N^R$ & 4470~(61)& 1319~(16) \\
Skewed &  $\hat{f}_S$&4253~(55) & 1343~(15) \\
& $\hat{f}_{S,N}$ & 4424~(61)& 1321~(16) \\
 & $\hat{f}_{S,N}^R$ & 4539~(62) & 1326~(16)  \\
\hline
~ & $\hat{f}$ & 4152~(59) & 1193~(16)\\
Kurtotic& $\hat{f}_N^R$ & 3882~(55)&~994~(13) \\
Unimodal&  $\hat{f}_S$& 4125~(59) & 1183~(16)\\
& $\hat{f}_{S,N}$ & 3860~(56)& 1007~(13) \\
 & $\hat{f}_{S,N}^R$ & 3882~(15) & ~994~(13)  \\
\hline
\end{tabular}

\end{center}

\newpage
\begin{center}

{\sc Table 1 continued}
\vspace{6pt}

\begin{tabular}{c|c|c|c}
\hline
 & & $n=100$ & $n=500$ \\
Density & Estimator & Mean Min.\ ISE  & 
Mean Min.\ ISE  \\
 & & (S.E.)  &
 (S.E.)  \\
\hline
~ & $\hat{f}$ & 4908~(110) & 1542~(32) \\
 & $\hat{f}_N^R$ & 2701~(~71)& ~737~(18) \\
Outlier&  $\hat{f}_S$& 4775~(102) &  1464~(30) \\
 & $\hat{f}_{S,N}$ & 3523~(~77)& ~940~(20) \\
 & $\hat{f}_{S,N}^R$ & 2926~(~70) & ~749~(18)  \\
\hline
~ & $\hat{f}$ & 717~(13)  & 223~(4)\\
& $\hat{f}_N^R$ & 658~(15)& 171~(4) \\
Bimodal & $\hat{f}_S$ &  702~(15) & 209~(4)\\
 & $\hat{f}_{S,N}$ & 625~(15)& 162~(4) \\
 & $\hat{f}_{S,N}^R$ & 639~(15) & 162~(4)  \\
\hline
~ & $\hat{f}$ & 1053~(19) & 313~(5)\\
Separated& $\hat{f}_N^R$ & ~711~(16)&  169~(4) \\
Bimodal &  $\hat{f}_S$&1021~(19) & 303~(5)\\
& $\hat{f}_{S,N}$ & ~813~(15)& 211~(4) \\
 & $\hat{f}_{S,N}^R$ & ~696~(16) & 165~(3)  \\
\hline
~ & $\hat{f}$ & 934~(15) & 299~(5) \\
Skewed& $\hat{f}_N^R$ & 924~(16)  & 269~(5)\\
Bimodal &  $\hat{f}_S$ & 950~(17) & 298~(5)\\
 & $\hat{f}_{S,N}$ & 952~(17)& 270~(5) \\
 & $\hat{f}_{S,N}^R$ & 970~(18) & 272~(5)  \\
\hline
%\end{tabular}
%\end{center}
%\newpage
%\begin{center}
%{Table 1 continued}
%\vspace{6pt}
%\begin{tabular}{c|c|c|c}
%\hline
% & &  $n=100$ & $n=500$ \\
%Density & Estimator & Mean Min.\ ISE  & 
%Mean Min.\ ISE  \\
% & & (S.E.) & 
% (S.E.)   \\
%\hline
~ & $\hat{f}$ & 864~(13) &  284~(4)\\
& $\hat{f}_N^R$ & 813~(13)&  268~(4) \\
Trimodal & $\hat{f}_S$&  852~(15)  & 280~(4)\\
 & $\hat{f}_{S,N}$ & 801~(14)& 269~(4) \\
 & $\hat{f}_{S,N}^R$ & 816~(15) & 271~(4)  \\
\hline
~ & $\hat{f}$ & 3652~(36) &1110~(11)  \\
 & $\hat{f}_N^R$ & 3754~(36)& ~994~(11) \\
Claw & $\hat{f}_S$ &  3666~(36) & 1108~(11) \\
 & $\hat{f}_{S,N}$ & 3693~(37)& ~995~(11) \\
 & $\hat{f}_{S,N}^R$ & 3782~(37) & ~994~(11)  \\
\hline
\end{tabular}

\end{center}

\vspace{.5em}
}

``Oracle'' versions of each of the
estimators under consideration were computed for each sample by empirical
calculation of the bandwidth that minimised the integrated squared error (ISE)
between estimate and true density. As a single global accuracy measure, we took
the ISE of the resulting (optimised) estimators averaged over simulations.
These (with standard errors in brackets) are given in Table 1.
Comparisons are, therefore, made in a ``best-case'' scenario, separating
estimator and bandwidth selection problems. It is 
possible that when data-based bandwidth selection is taken into account,
varying degrees of difficulty in choice of bandwidth {\it could} change
the relative merits of the procedures (but the basic method and 
bandwidth selection problems would then be conflated). One would not expect,
however, that a bandwidth selector for $\hat{f}_{S,N}^R$, say, would be 
so much more effective than one for $\hat{f}_N$, say, that the 
improvement of the former over the latter would be greatly enhanced in 
that case. Indeed, results of, for example, Jones (1992) suggest that the 
quality of data-based bandwidth selectors will in fact go down with 
improved performance of the basic estimator.

Five estimators are compared in the table. Each of those 
with a parametric component employs 
the normal distribution in that guise. Results for the basic kernel 
density
estimator (1.3) and the renormalised Jones, Linton and Nielsen (1995) estimator
are the same as in Jones and Signorini (1997). We do not bother with the raw
estimator $\hat{f}_N$ because previous authors have shown the renormalisation
to be uniformly advantageous. It turns out that although $\int \hat{f}_S(x)dx
\neq 1$, the difference from unity is very small and the effect on performance
of renormalising $\hat{f}_S$ in practice is almost negligible. For this reason,
we exhibit only one version of the Hjort and Glad estimator in Table 1, and this
is the raw form $\hat{f}_S$. Renormalisation does make a noticeable difference
in the case of $\hat{f}_{S,N}$, however. Results for
both raw and renormalised HOBSKDEs are given since, unexpectedly, neither
$\hat{f}_{S,N}$ nor $\hat{f}_{S,N}^R$ dominates the other. If one of the two
has to be preferred, perhaps renormalisation wins because, while the difference
between the two is small in many cases, where there are the most substantial
differences (Outlier, Separated Bimodal) the renormalised version is the better
of the two.

Consider the $n=100$ results.
For seven of the ten densities, the basic kernel estimator $\hat{f}$ is
improved upon by all three alternative estimators. Away from the Gaussian
distribution, the degree of improvement made by the Hjort and Glad estimator
is often small. Then again, when $\hat{f}$ does relatively well (Strongly
Skewed, Skewed Bimodal, Claw), so does $\hat{f}_S$. (More on performance
of SKDEs will be found in Hjort, Jones and Storvik.) At the Gaussian density, 
there are great
improvements over $\hat{f}$ from all alternative estimators, but interestingly
it is $\hat{f}_N^R$ and $\hat{f}_S$ that lead the way, with neither (2.1) nor
(2.2) able to take quite as much advantage of the situation. Elsewhere,
$\hat{f}_N^R$ and $\hat{f}_{S,N}^R$ have broadly comparable
performance (generally better than that of $\hat{f}_S$) with overall a slight
preference for $\hat{f}_N^R$.

Similar relative performances can be observed for the $n=500$ sample size
although (i) renormalisation of $\hat{f}_{S,N}$ is now uniformly no worse than
the raw version, and (ii) performance of $\hat{f}_{S,N}^R$ does generally get a
little better, in line with its asymptotic justification.

\vspace{1em}
\begin{center}
{\bf 4. Conclusions}
\end{center}
\vspace{1em}

The new estimators $\hat{f}_{S,N}$ or $\hat{f}_{S,N}^R$ realise {\it some}
of their theoretical potential in practice. For the smaller sample size
($n=100$), they perform well
relative to $\hat{f}$ and, in general, relative to $\hat{f}_S$, although they
do not perform so well as the latter at the parametric model. Performance is,
however, somewhat disappointing in that they are unable to
improve in general on $\hat{f}_N^R$. Notice, however, that the
good performance of $\hat{f}_N^R$ near the normal model would not 
transfer to other parametric models which may be used as targets; the 
semiparametric estimators would transfer good performance readily to
alternative vehicle models. For the larger sample size 
$(n=500$), $\hat{f}_{S,N}^R$
comes more into its own, although still it is unclear whether it would be
practically worthwhile to prefer $\hat{f}_{S,N}^R$ to $\hat{f}_N^R$ 
in this case.
Caveats concerning automatic bandwidth selection remain in place.

The
proposals of this paper are not unique in their theoretical properties,
but it is not clear that it is worth investigating further methods in the 
same class of HOBSKDEs.

\vspace{1em}
\begin{center}
{\bf References}
\end{center}
\vspace{1em}

{\footnotesize

\begin{description}
\item {\sc Hjort, N.L. and Glad, I.K.} (1995). Nonparametric density estimation 
with a parametric start. {\it Ann. Statist.} {\bf 23}, 882-904.
\vspace{-.5em}
\item {\sc Jones, M.C.} (1992). Potential for automatic bandwidth choice in 
variations on kernel density estimation. {\it Statist. Probab. Lett.} {\bf 
13}, 351-356. 
\vspace{-.5em}
\item {\sc Jones, M.C., Linton, O. and Nielsen, J.P.} (1995). A simple and 
effective bias reduction
method for kernel density estimation. {\it Biometrika} {\bf 82}, 327-338.
\vspace{-.5em}
\item {\sc Jones, M.C. and Signorini, D.F.} (1997). A comparison of higher order
bias kernel density estimators. {\it J. Amer. 
Statist. Assoc.} {\bf 92}, 1063-1073.
\vspace{-.5em}
\item {\sc Marron, J.S. and Wand, M.P.} (1992). Exact mean integrated squared 
error. {\it Ann. Statist.} {\bf 20}, 712-736.
\vspace{-.5em}
\item {\sc Scott, D.W.} (1992). {\it Multivariate Density Estimation: Theory, 
Practice, and Visualization}. Wiley, New York.
\vspace{-.5em}
\item {\sc Wand, M.P. and Jones, M.C.} (1995). {\it Kernel Smoothing}.
Chapman and Hall, London.
\end{description}

\vspace{2em}

\noindent
\begin{tabular}[t]{l}
{\sc M.C. Jones}\\
{\sc Department of Statistics}\\
{\sc The Open University}\\
{\sc Walton Hall}\\
{\sc Milton Keynes MK7 6AA}\\
{\sc U.K.}\\
email: m.c.jones@open.ac.uk\\
\end{tabular}
\quad \quad 
\begin{tabular}[t]{l}
{\sc D.F. Signorini} \\
{\sc Biostatistics}\\
{\sc Quintiles Scotland Ltd}\\
{\sc Bathgate}\\
{\sc West Lothian EH48 2EH}\\
{\sc U.K.}\\
email: dsignori@qedi.quintiles.com\\
\end{tabular}
\quad \quad 
\begin{tabular}[t]{l}
{\sc N.L. Hjort}\\
{\sc Department of Mathematics}\\
{\sc University of Oslo}\\
{\sc P.B. 1053 Blindern}\\
{\sc N--0316 Oslo 3}\\
{\sc Norway}\\
email: nils@math.uio.no\\
\end{tabular}

}

\end{document}